\documentclass[11pt]{article}

\usepackage{amsmath, amssymb, amsthm}
\usepackage[english]{babel}
\usepackage[T1]{fontenc}
\usepackage{enumerate}
\numberwithin{equation}{section}
\newcommand{\be}{\begin{equation}}
\newcommand{\ee}{\end{equation}}
\newcommand{\bd}{\begin{displaymath}}
\newcommand{\ed}{\end{displaymath}}
\newcommand{\ba}{\begin{eqnarray}}
\newcommand{\ea}{\end{eqnarray}}

\newenvironment{equations}{\equation\aligned}{\endaligned\endequation}

\def\R{{I \!\! R}}

\def\v12{(v-w)}

\def\({\left(}
\def\){\right)}

\def\bgr#1\egr{{\allowdisplaybreaks\begin{gather}#1\end{gather}}}
\def\bma#1\ema{{\allowdisplaybreaks\begin{align}#1\end{align}}}

\def\oplem#1{\begin{lemma}\, {\rm #1}\, \it }
\def\cllem{\end{lemma}\rm \par }
\def\opthm#1{\begin{theorem}\, {\rm #1}\, \it }
\def\clthm{\end{theorem}\rm \par }

\def\R{{\rm{I}\! \rm{R}}}

\def\pRR{\hbox{{\tiny \rm I}\kern-.1em\hbox{{\tiny \rm R}}}}

\def\NN{\hbox{I\kern-.2em\hbox{N}}}

\newcommand{\mR}{\mathcal H}
\newcommand{\mN}{\mathcal N}
\newcommand{\mH}{\mathcal H}

\newcommand{\mI}{\mathcal I}

\def\R{\mathbb{R}}

\newcommand{\fer}[1]{(\ref{#1})}
\newcommand{\bq}{\begin{equation}}
\newcommand{\eq}{\end{equation}}
\def\bqa{\begin{eqnarray}}
\def\eqa{\end{eqnarray}}
\def\bd{\begin{displaymath}}
\def\ed{\end{displaymath}}

\theoremstyle{remark}

\theoremstyle{definition}

\begin{document}

\title{The information-theoretic meaning of Gagliardo--Nirenberg type inequalities}
\date{}

 \author{Giuseppe Toscani \thanks{Department of Mathematics, University of Pavia, and IMATI CNR, via Ferrata 1, 27100 Pavia, Italy.
\texttt{giuseppe.toscani@unipv.it} }}

\maketitle

\begin{flushleft}
\emph{Dedicated to the memory of Emilio Gagliardo}
\end{flushleft}

\vskip 2cm
\begin{center}\small
\parbox{0.85\textwidth}{
\textbf{Abstract.} Gagliardo--Nirenberg inequalities are interpolation inequalities which were proved independently by Gagliardo  and Nirenberg  in the late fifties. In recent years, their connections with theoretic aspects of information theory and nonlinear diffusion equations allowed to obtain some of them in optimal form, by recovering both the sharp constants and the explicit form of the optimizers.  In this note, at the light of these recent researches, we review the main connections
between Shannon-type entropies,  diffusion equations and a class of these inequalities.
\medskip

\textbf{Keywords.} Heat equation, nonlinear diffusions, Shannon entropy, R\'enyi entropy, Fisher-type informations, Gagliardo--Nirenberg type inequalities.}
\end{center}

\medskip

\section{Introduction}
\label{intro}
Gagliardo--Nirenberg inequalities, proved independently by Gagliardo in \cite{Gagl} and Nirenberg \cite{Nire} in the late fifties, are interpolation inequalities of the form
 \be\label{GN1}
 \| u\|_{L^p(\R^n)} \le K_{GN} \|\nabla u\|_{L^2(\R^n)}^\theta \| u\|_{L^q(\R^n)} ^{1-\theta},
  \ee
for $u \in W^{1,q}(\R^n)$.  In \fer{GN1} $K_{GN}$ is a positive constant, $n>2$ and $1<q<p<2^* = 2n/(n-2)$, or $n=1,2$ and $1<q<p$. Last $\theta = [2n(1-q/p)]/[2n-q(n-2)]$. Here we denoted
 \[
W^{1,q}(\R^n) = \{u \in L^q(\R^n): \nabla u \in L^2(\R^n)\}.
\]
and the $L^2$-norm of $\nabla u$ has been considered for simplicity, though in general, when $n >2$ the Gagliardo--Nirenberg inequalities can be stated with the $L^r$-norm of $\nabla u$, where $1 < r < n$. 

The problem of finding the sharp constants and optimal functions for these inequalities has attracted many researchers in the past years \cite{Ag1,Agu,CNV,DD,Lev}. Also, results on their stability with respect to optimizers has been  studied in various aspects (cf. \cite{CF,DT3} and the references therein). 

A maybe not so well-known aspect of a special class of these inequalities,  is their relationships  with information theory and entropy-type inequalities \cite{ST,Tos4,Tos5}. 

In information theory, the main examples are the celebrated Shannon entropy power inequality, formulated by Shannon in his pioneering paper \cite{Sha}, and the Blachman--Stam inequality \cite{Bla,Sta}.
Given a random vector $X$ in $\R^n$, $n \ge 1$ with density $f(x)$, let
 \be\label{shan}
\mR(X) = \mR(f) = - \int_{\R^n} f(x) \log f(x)\, dx
 \ee
denote its entropy functional (or Shannon entropy). Together with entropy, Shannon  introduced the concept of  entropy power, defined by 
  \be\label{ep}
 \mN(X) = \mN(f) = \exp\left(\frac 2n \mR(X)\right).
 \ee
The entropy power is built to be in a suitable sense \emph{linear} at  Gaussian random vectors. To be precise, let $Z_\sigma = N(0,\sigma I_n)$ denote the $n$-dimensional Gaussian random vector having mean vector $0$ and covariance matrix $\sigma I_n$, where $I_n$ is the identity matrix. Then, the entropy power of $Z_\sigma$ is linear in terms of the variance, since $\mN(Z_\sigma) = \sigma$. Shannon entropy power inequality, first rigorously proven by Stam \cite{Sta} (cf. also
\cite{Cos, GSV, GSV2, Rio, Tos4} for other proofs and extensions) gives a lower bound on
Shannon entropy power of the sum of independent random variables $X, Y$ in $\R^n$ with densities
 \be\label{entr}
\mN(X+Y) \ge \mN(X) + \mN(Y),
 \ee
with equality if and only if $X$ and $Y$ are Gaussian random vectors with proportional
covariance matrices. In other words, Shannon entropy power characterizes an extremal property of Gaussian functions with respect to convolutions.

Likewise, Blachman--Stam inequality is concerned with the behavior of the Fisher information
with respect to convolutions. Historically, it was the key argument to prove
Shannon entropy power inequality  \cite{Bla,Sta}. 
Given the $n$-dimensional random vector $X$ of  probability density $f(x)$, let
 \be\label{fis}
\mI(X) = \mI(f) = \int_{\pRR^n} f(x)\left|\nabla\log f(x)\right|^2\, dx = \int_{\{f>0\}} \frac{|\nabla f(x)|^2}{f(x)} \, dx.
 \ee
define its Fisher measure of information. 
Analogously to the entropy power inequality \fer{entr}, Blachman--Stam inequality furnishes a lower bound on the Fisher information of the sum of independent random variables $X, Y$ in $\R^n$ with densities
 \be\label{BS}
\frac 1{\mI(X+Y)} \ge \frac 1{\mI(X)} + \frac 1{\mI(Y)},
 \ee
with equality if and only if $X$ and $Y$ are Gaussian random vectors with proportional
covariance matrices. Hence, analogously to \fer{entr}, inequality \fer{BS} characterizes a further extremal property of Gaussian functions with respect to convolutions. Note that, similarly to the entropy power, $1/{\mI(X)}$ shares the same \emph{linearity property} of the entropy power. Indeed  $1/{\mI(X)}= n\sigma$ when $X=Z_\sigma$.

In Shannon's theory of information, the proof of these inequalities  has benefited from the close links between entropies and diffusion equations. These links are perfectly understood in the linear case, where, in addition to its paramount importance in physical applications, the linear heat equation is known to represent a profitable tool to obtain mathematical inequalities in sharp form \cite{B2,BB,BBC,B1,Tos3,Tos4}. 

This original way of using heat equation dates back to the years between the late fifties to mid sixties, exactly at the same time in which Gagliardo \cite{Gagl} and Nirenberg \cite{Nire} proved the interpolation inequalities that bear their name. The pioneering application of heat equation to the finding of analytic inequalities is due to Stam \cite{Sta} in 1959 (cf. also Blachman \cite{Bla}), who used the link between entropy and Fisher measure of information, obtained by deriving Shannon entropy along the solution to the heat equation, to find a rigorous proof of inequality \fer{entr}. It is interesting to remark that the same link  was used independently by Linnik \cite{Lin} (his research appeared in the same year of Stam's work \cite{Sta}) to get a new information-theoretic proof of the central limit theorem of probability theory. 

Also, some years later, the heat equation has been used in the context of kinetic theory of rarefied gases by McKean \cite{McK} to investigate the large-time behavior of Kac caricature of a Maxwell gas. There, various monotonicity properties of the derivatives of Shannon entropy along the solution to the heat equation have been obtained.

These contributions made popular with the information community the deep link between Shannon entropy \fer{shan}
and the solution of the heat equation  posed in the whole space $\R^n$
  \be\label{heat}
\frac{\partial u(x,t)}{\partial t} =  \Delta u(x,t),
 \ee
as soon as the initial datum is assumed to be a probability density.

More recently \cite{DT2,ST,Tos5}, it has noticed that a similar link can be established between  R\'enyi entropy of order $p$ and the nonlinear diffusion of order $p$,  posed in the whole space $\R^n$
\begin{equation}
{\partial v(x,t) \over \partial t} =   \Delta v^p(x,t),
 \label{poro}
\end{equation}
still with the initial datum assumed to be a probability density, in the range $p>(n-2)/n$.

Given a random vector $X$ in $\R^n$, $n \ge 1$ with density $f(x)$, and a positive constant $p$,
R\'enyi entropy of order $p$ of  $X$  is defined by \cite{DCT}
 \be\label{re}
 \mR_p(X)=
\mR_p(f) = \frac 1{1-p} \log\left( \int_{\pRR^n} f^p(y) \, dy \right).
 \ee
This concept of entropy was introduced by R\'enyi in \cite{Ren} for a discrete pro\-ba\-bility measure to generalize the classical logarithmic entropy, by maintaining at the same time most of its properties. Indeed,  R\'enyi entropy of order $1$,  defined as the limit as $p \to 1$ of $\mR_p(f)$ coincides with Shannon entropy.
Therefore, the standard (Shannon) entropy  of a probability density \cite{Sha} is included in the set of R\'enyi entropies, and it is identified with R\'enyi entropy of index $p = 1$. 

The deep link between nonlinear diffusions and R\'enyi entropies gave a new light to a certain class of inequalities already present in the literature.
Indeed, inequalities for R\'enyi entropies and generalized Fisher information measures were considered before \cite{LYZ,LYZ2}, without resorting to any connection with nonlinear diffusions.

While the derivation of sharp inequalities involving entropies, and the (eventual) characterization of the extremal densities, represents one of the main objectives of information theory \cite{CT,DCT}, often this point of view and the subsequent results do not spread automatically to other mathematical communities. 

The aim of this note is to outline the close relationships between entropy type and theoretic inequalities. In particular, we will show that a certain class of Gagliardo--Nirenberg type interpolation inequalities can be viewed as extremal properties of R\'enyi measures of information. Consequently, their derivation can be obtained by resorting to arguments which are strongly based on their information-theoretic meaning.  These results are rooted in many contributions from the field of information theory \cite{Cos,CT,CHV,DCT,LYZ,LYZ2,Vil}, mass transportation  and nonlinear diffusion equations \cite{BDV,BDGV,CaTo,CNV,CHV,DD}. Also, details on some of these results can be extracted from various recent papers of this author \cite{CTo,DT1,DT2,ST,Tos2,Tos3,Tos4,Tos5}. 

\section{Heat equation and related entropies}
\label{sec:heat}
The goal of this Section is to enlighten, by means of a simple example, how the solution to the linear heat equation can play a role in recovering information-type inequalities in sharp form. To start with, let us briefly recall some elementary but basic result about the linear diffusion equation. Further details can be found in any introductory book on partial differential equations (cf. for example the classical treatise of Evans \cite{Eva}). 

Let $u(x,t)$ denote the solution to the initial value problem of equation \fer{heat}, corresponding to an initial value $u_0(x)$ that is a probability density with finite second moment. The solution to the initial value problem for equation \fer{heat} is given by the convolution product of the initial density $u_0(x)$ with the Gaussian density $M_{2t}$, being  
 \be\label{gau}
 M_{\sigma}(x) = \frac 1{(2\pi \sigma)^{n/2}} \exp\left\{-\frac{x^2}{2\sigma} \right\}.
 \ee
the Gaussian density  of the $n$-dimensional Gaussian random vector $Z_{\sigma}$ of mean vector zero and covariance matrix  $\sigma I$ (with variance $n\sigma$).

It is well-known (and in any case immediate to verify) that the solution to the heat equation is such that mass and momentum are preserved in time, so that
 \[
 \rho(u(t)) = \int_{\pRR^n} u(x,t) \, dx =  \int_{\pRR^n} u_0(x) \, dx = 1, 
  \]
  and
  \[
  m(u(t)) = \int_{\pRR^n} x u(x,t) \, dx =  \int_{\pRR^n}x u_0(x) \, dx = m(u_0).
 \]
Differently, the second moment
\[
E(u(t)) =  \int_{\pRR^n} |x|^2 u(x,t) \, dx
\] 
is growing linearly in time, at a rate given by
 \be\label{m2}
 \frac d{dt}E(u(t)) = 2n\rho(u(t))= 2n.
 \ee
Note that the growth of the second moment depends on the initial value only through its mass density. This implies that, in correspondence to each initial datum of unit mass, independently of its shape, the rate of growth has the same value $2n$. However, a more refined estimate can be obtained by resorting to the forthcoming argument, which takes advantage of a particular rewriting of the Laplace operator. 

Since the initial value $u_0(x)$ is a probability density, $u_0(x) \ge 0$, and the maximum principle insures that the solution $u(x,t)$ remains nonnegative for all subsequent times $t >0$. Thus, we can write the diffusion equation \fer{heat} in the alternative form
 \be\label{heat2}
 \frac{\partial u(x,t)}{\partial t} = \sum_{i=1}^n \frac{\partial}{\partial x_i}\left( u(x,t) \frac{\partial}{\partial x_i}\log u(x,t)\right).
 \ee
 For any given time $t >0$, integration by parts gives
 \begin{equations}\label{inte}
 \frac d{dt} E(u(t)) =& \frac d{dt} \int_{\pRR^n} |x|^2 u(x,t)\, dx = \int_{\pRR^n} |x|^2  \frac{\partial u(x,t)}{\partial t}\, dx = \cr
 & \int_{\pRR^n} |x|^2\sum_{i=1}^n  \frac{\partial}{\partial x_i}\left( u(x,t) \frac{\partial}{\partial x_i}\log u(x,t)\right)\, dx = \cr 
& -2 \int_{\pRR^n}\sum_{i=1}^n x_i u(x,t) \frac{\partial}{\partial x_i}\log u(x,t)\, dx.
 \end{equations}
By Cauchy--Schwarz inequality 
 \be\label{CS1}
 \left| \int_{\pRR^n}\sum_{i=1}^n x_i \frac{\partial}{\partial x_i}\log u(x,t)\, u(x,t) \, dx\right| \le E(u(t))^{1/2} \mI(u(t))^{1/2},
 \ee
being $\mI(u(t))$ the Fisher information of the solution (cf. definition \fer{fis}).
We remark that, as proven in \cite{LT}, Lemma $\rm 2.1$, $\mI(u(t))$  is bounded as soon as $t >0$. 
Inequality  \fer{CS1} then implies that the growth of the square root of the second moment of the solution to the heat equation can not exceed the square root of the Fisher information of the solution itself. Indeed, substituting the right-hand side of \fer{CS1} into the last line of \fer{inte} we get
\be\label{gro1}
  \frac d{dt} \sqrt{E(u(t))} \le \sqrt{\mI(u(t))}.
 \ee
Moreover, using \fer{m2} into \fer{gro1} we conclude with the information-type inequality \cite{DCT}
 \be\label{fis1}
 E(u(t))\mI(u(t)) \ge n^2. 
 \ee
Note that both inequality \fer{CS1} and the differential inequality \fer{gro1} become equalities when evaluated in correspondence to the Gaussian $M_{2t}$, $t >0$, fundamental solution of the heat equation \fer{heat}.  

Indeed, for the Gaussian function $M_{\sigma}(x)$, it holds
 \[
 -\int_{\pRR^n}\sum_{i=1}^n x_i M_{\sigma}(x) \frac{\partial}{\partial x_i}\log M_{\sigma}(x)\, dx = n,
 \] 
while
 \[
 \mI(M_{\sigma}) = \frac n\sigma.
 \]
Let us draw some conclusions from the previous example. We studied the evolution of a time dependent \emph{functional} of the solution  to the heat equation (the second moment) which grows  linearly with respect to time. In the case of the linear diffusion, this linear growth is common to all solutions. 

Since the growth is linear, the time derivative of the functional is constant. In particular, this is true if we consider the evolution of the second moment of the fundamental solution $M_{2t}$. Note however that time appears in the fundamental solution as a dilation parameter, where, as usual, for any given density $f(x)$ and positive constant $a$,  we define the \emph{dilation} of $f$ by $a$, as the mass-preserving scaling
 \be\label{scal}
 f(x) \to f_a(x) = {a^n} f\left( {a} x \right).
 \ee
Consequently, the time derivative of the functional, which is constant in time, is invariant with respect to dilations (scale invariant). Hence, we can take the limit as $t \to 0^+$ of the left--hand side of inequality \fer{fis1} to obtain, for any probability density $f(x)$, $x \in \R^n$, with bounded Fisher information the scale invariant inequality
 \be\label{fis2}
 E(f)\mI(f) \ge n^2.
 \ee
Note that the constant $n^2$ in inequality \fer{fis2} is sharp on the set of probability densities. This value is reached in correspondence to Gaussian densities, which are the unique minimizers. Moreover, for any given initial value $u_0(x)$ which is a probability density, the Gaussian density is the unique initial value for which equality in \fer{gro1} holds. 

A further example will help to understand the remarkable potential of this idea. Let us study the time evolution of Shannon entropy, as defined in \fer{shan}, of the solution to the heat equation. Using the heat equation written in the form \fer{heat2} we easily conclude with DeBruijn's identity
\be\label{e-s}
 \frac {d}{dt} \mR(u(t)) =  \mI(u(t)), \quad t >0,
 \ee
that links Shannon entropy with Fisher measure of information via the heat equation.
Now, consider that Shannon entropy of the fundamental solution to the heat equation  has the value
 \be\label{hm}
 \mR(M_{2t}) = \frac n2\log(4\pi e t).
 \ee
Hence, it grows logarithmically with time. To restore the linear growth we need to consider the quantity
 $ \exp\left\{ 2 \mR(f)/n\right\}$,  namely Shannon entropy power of a pro\-ba\-bility density defined in \fer{ep}. This clarifies (from the point of view of the linear diffusion equation) why the entropy power is the right quantity to study to get inequalities. 
 
As in the previous example,  let us compute the subsequent derivatives of Shannon entropy power.  These computations lead to a well-known result in information theory, known as the \emph{concavity of entropy
power}. If $u(x,t)$ is a solution to the heat
equation \fer{heat}, corresponding to an initial datum $u_0(x)$ that
is a probability density, then its entropy power is a concave function of time
 \be\label{concav}
\frac{d^2}{dt^2}\mN(u(t)) \le 0.
 \ee
Moreover, equality in \fer{concav} holds if and only if $u(x,t)$
coincides with the Gaussian density of variance $2t$, namely the
fundamental solution to the heat equation.  Inequality \fer{concav}
has been first obtained by Costa \cite{Cos}. Simplified proofs were subsequently done by Dembo \cite{Dem} by means of a direct application of Blachman--Stam inequality \fer{BS}, and Villani \cite{Vil}, who made use of an argument 
introduced by McKean \cite{McK} in his paper on Kac caricature of a Maxwell gas.

Maybe the most important consequence of \fer{concav} is the sharp inequality which can be extracted from the first derivative of the entropy power \cite{Tos2}.  Indeed, analogously to the previous example, the first derivative is both scale invariant, and nonincreasing in consequence of the concavity property. Since
 \be\label{pd}
 \frac{d}{dt}\mN(u(t)) = \frac 2n\mN(u(t))\mI(u(t)), 
 \ee
 invariance under dilation coupled with concavity allow to reach the lower bound by taking the limit as $t \to \infty$ of the right-hand side of equation \fer{pd}. This implies the so-called \emph{isoperimetric inequality for entropies} \cite{DCT, Tos2}. 
For any probability density $f(x)$ in $\R^n$, this scale invariant
inequality asserts that
  \be\label{b51}
\mN(f) \, \mI(f) \ge  2\pi e n.
 \ee
As for \fer{fis1}, the constant $2\pi e n$ is sharp on the set of probability densities. This value is reached in correspondence to Gaussian densities, which are the unique minimizers. Moreover, for any given initial value $u_0(x)$ which is a probability density, the Gaussian density is the unique initial value for which equality in \fer{concav} holds. 

It is interesting to remark that from the isoperimetric inequality for entropies \fer{b51} one can obtain, among other consequences, the logarithmic Sobolev inequality in scale invariant form, and Nash's inequality \cite{Nash} with an asymptotically sharp constant \cite{Tos2}. In particular, Nash's interpolation inequality reads 
\be\label{nass1}
\left( \int_{\R^n} g^2(x) \, dx \right)^{1+ 2/n} \le   \frac 2{\pi e
n} \left( \int_{\R^n}|g(x)\, dv \right)^{4/n} \int_{\R^n} |\nabla
g(x) |^2 \, dx.
 \ee

\section{From R\'enyi's entropies to Gagliardo--Nirenberg inequalities}
\label{sec2}

The arguments of Section \ref{sec:heat}, enlightening the role of the linear heat equation in recovering sharp inequalities, can be generalized to the nonlinear diffusion equations \fer{poro} and relative entropies. 

Similarly to Section \ref{sec:heat}, our analysis will be restricted to initial data which are probability densities with finite variance. This allows to include both the case $p>1$, usually known with the name of porous medium equation, and a limited range of exponents when $p<1$, the fast diffusion equation. In dimension $n\ge 1$, the second moment of the solution to equation \fer{poro} is  bounded if $p > \bar p
$ with $\bar p = n/(n+2)$.  The particular subinterval of $p$ is
motivated by the existence of a precise solution, found by
Zel'dovich, Kompaneets and Barenblatt in the fifties (briefly
called here Barenblatt solution) \cite{Ba1,Bar,ZK}, which serves
as a model for the asymptotic behavior of a wide class of solutions
with finite second moment. In the case $p>1$ (see \cite{BDV} for $p
<1$) the Barenblatt (also called self-similar or generalized
Gaussian solution)  departing from $x=0$ takes the self-similar form
 \be\label{ba-self}
 M_p(x,t) := \frac 1{t^{n/\mu}} \tilde M_p\left (\frac x{t^{1/\mu}}\right),
 \ee
 where
\[
\mu =2+n(p-1)
\]
and  $\tilde M_p(x)$ is the time-independent function
\begin{equation}
 \tilde M_p(x)=\big(C - \lambda\, |x|^2   \big)^{\frac
1{p-1}}_+ . \label{ba}
 \end{equation}
In \fer{ba} $(s)_+ = \max\{s, 0\}$,
$\lambda = \frac 1{2\mu}\,\frac{p-1}{p}$, and the constant $C$  is chosen to fix the mass of the source-type Barenblatt solution equal to one.

Since we are interested in presenting the information-theoretic aspect of inequalities, in what follows, we will often proceed formally. However, at the price of an increasing number of technical details, the results can be fully justified on the basis of classical results of existence, uniqueness and regularity of the solution to the initial value problem  \cite{BDV, ST,Vaz}. 
In analogy with the heat equation, the solution to the nonlinear diffusion \fer{poro} is such that mass and momentum are preserved in time
 \[
 \rho(v(t)) = \int_{\pRR^n} v(x,t) \, dx =  \int_{\pRR^n} v_0(x) \, dx = 1, 
  \]
  and
  \[
  m(v(t)) = \int_{\pRR^n} x v(x,t) \, dx =  \int_{\pRR^n}x v_0(x) \, dx = m(v_0).
 \]
Moreover, while the second moment  $E(v(t))$ of equation \fer{poro} increases in time from $E_0 = E(v_0)$, differently from the heat equation its
evolution is now given by the nonlinear law
 \be\label{ene}
\frac{dE(v(t))}{dt} = 2n \int_{\pRR^n}v^p(x,t)\, dx \ge 0,
 \ee
which, since $p\not=1$, it is not explicitly computable.

Note that in the nonlinear case the growth of the second moment varies with time, and depends on the solution itself.  Nevertheless, one can try to apply to equation \fer{poro} the same strategy applied to the linear heat equation. 
Similarly to the case $p = 1$, we can write the nonlinear diffusion equation in the alternative form
 \be\label{heat22}
 \frac{\partial v(x,t)}{\partial t} = \sum_{i=1}^n \frac{\partial}{\partial x_i}\left(p\, v(x,t) \frac{\partial}{\partial x_i}\frac{v^{p-1}(x,t)}{p-1}(x,t)\right).
 \ee
 For any given time $t >0$, integration by parts gives
 \begin{equations}\label{inte2}
 \frac d{dt} E(v(t)) =& \frac d{dt} \int_{\pRR^n} |x|^2 v(x,t)\, dx = \int_{\pRR^n} |x|^2  \frac{\partial v(x,t)}{\partial t}\, dx = \cr
 & \int_{\pRR^n} |x|^2\sum_{i=1}^n  \frac{\partial}{\partial x_i}\left(p\, v(x,t) \frac{\partial}{\partial x_i}\frac{v^{p-1}(x,t)}{p-1}(x,t)\right)\, dx = \cr 
& -2\,p\, \int_{\pRR^n}\sum_{i=1}^n x_i v(x,t) \, \frac{\partial}{\partial x_i} \frac{v^{p-1}(x,t)}{p-1}(x,t)\, dx.
\end{equations}
By Cauchy--Schwarz inequality 
 \begin{equations}\label{CS2}
 &\left| \int_{\pRR^n}\sum_{i=1}^n x_i \frac{\partial}{\partial x_i}\, v(x,t)\frac{v^{p-1}(x,t)}{p-1}(x,t) \, dx\right| \le \\
& E(v(t))^{1/2} \left(\int_{\R^n} v(x,t) \left|\frac{\nabla v^{p-1}(x,t)}{p-1} \right|^2 \, dx\right)^{1/2}= \\
& \frac 1p E(v(t))^{1/2}\left(\int_{\R^n}\frac{\left|\nabla v^p(x,t)\right|^2}{v(x,t)} \, dx\right)^{1/2}. 
 \end{equations}
For a given $n$-dimensional random vector $X$ of probability density $f(x)$, let us define the Fisher information of order $p$ as \cite{ST}
 \be\label{fis-r}
 \mI_p(X) = \mI_p(f) := \frac 1{\int_{\pRR^n} f^p(x) \, d x}
   \int_{\{f>0\}} \frac{|\nabla f^p(x)|^2}{f(x)} \, d x.
 \ee 
We remark that, as $p \to 1$ the Fisher information of order $p$ tends towards the classical Fisher information defined by \fer{fis}. 
Then, from equality \fer{ene} we obtain the bound
 \be\label{fis3}
 E(v(t))I_p(v(t)) \ge n^2\, \int_{\pRR^n} v^p(x,t) \, d x.
 \ee
By taking the limit as $t \to 0^+$ of both sides of inequality \fer{fis3}, it follows that any probability density $f(x)$, $x \in \R^n$, with bounded variance satisfies the inequality
\be\label{fis33}
 E(f)I_p(f) \ge  n^2\, \int_{\pRR^n} f^p(x) \, d x. 
 \ee
As before, equality is reached in \fer{fis33} if and only if $f= \tilde M_p$, the Barenblatt profile of order $p$. Note that, at difference with the linear case, when $p\not=1$, the derivative of the second moment can be directly related to R\'enyi entropy. Indeed, by \fer{re} 
 \be\label{der1}
  \int_{\pRR^n} f^p(x) \, d x = \exp\left\{ (1-p)\mH_p(f)\right\}.
   \ee 
Hence, we can rewrite inequality \fer{fis3} as an information-type inequality
\be\label{fis4}
 E(f)I_p(f) \ge  n^2\, \exp\left\{ (1-p)\mH_p(f)\right\}, 
 \ee
which contains inequality \fer{fis1} as $p \to 1$.

 The definition of the $p$-Fisher information is motivated by the following relationship.
Let us consider the evolution in time of R\'enyi entropy
of order $p$, defined in \fer{re}, along the solution of the nonlinear diffusion equation
\fer{poro}. Integration by parts then yields, for $t >0$
 \be\label{e-r}
 \frac {d}{dt} \mR _p(v(\cdot,t)) =  \mI_p(v(\cdot,t)).
 \ee 
When $p \to 1$, identity \fer{e-r} reduces to DeBruijn's identity \fer{e-s}
which connects Shannon entropy functional with the Fisher
information.

Before to proceed, let us recall some interesting information-theoretic consequences of inequality \fer{fis3} (alternatively \fer{fis4}). 
Among other properties,  R\'enyi entropy \fer{re} behaves as Shannon entropy \fer{shan} with respect to dilations of the probability density.
Indeed, for any $p \ge 0$ it holds
 \be\label{dil}
 \mR_p(f_a) = \mR_p(f) - n \log a.
 \ee
This property implies that, for any given probability density function $f$ in $\R^n$ such that both the second moment and  R\'enyi entropy are bounded, the functional
 \be\label{ok1}
 \Lambda(f) = \mH_p(f) - \frac n2 \log E(f)
 \ee
is invariant with respect to dilations.  Let us consider the time variation of  $\Lambda(v(t))$, where $v(x,t)$ is the solution to the nonlinear diffusion \fer{poro} corresponding to an initial value $u_0(x)$ which is a probability density with finite second moment. Clearly, thanks to inequality \fer{fis3} it holds
  \be\label{der6}
  \frac d{dt} \Lambda(v(t)) =  \mI_p(v(t)) - n^2 \frac{ \int_{\pRR^n} v^p(x,t) \, d x}{ E(v(t))} \ge 0.
  \ee
 Thus, $\Lambda(v(t))$ is increasing in time, which implies, grace to dilation invariance, that it will converge, as time $t\to \infty$ to the value obtained in correspondence to the Barenblatt profile \cite{Tos5}.
Hence, given  a probability density function $f(x)$, $x \in\R^n$, such that both the second moment and  R\'enyi entropy are bounded, it holds
\be\label{in9}
\mH_p(f) - \frac n2 \log E(f) \le  \mH_p(\tilde M) - \frac n2 \log E(\tilde M_p).
\ee
The same inequality has been proven in \cite{CHV,LYZ,LYZ2} by different methods.

Let
$\mN_p(f)$ denote the entropy power of $f$ associated to R\'enyi's
entropy of order $p$
 \be\label{epo}
\mN_p(f) = \exp\left\{ \left(\frac 2n + p-1\right)\mR_p(f)\right\}.
 \ee
Then, if  $p> n/(n+2)$, we rewrite \fer{in9} as the information type inequality
 \be\label{ep1}
 \frac{\mN_p(f)}{E(f)^{1+ n(p-1)/2}} \le
 \frac{\mN_p(B_{p,\sigma})}{E(B_{p,\sigma})^{1+ n(p-1)/2}}.
 \ee
We remark that the definition \fer{epo} of $p$-R\'enyi entropy power,
proposed recently in \cite{ST}, coincides with the classical
definition of Shannon entropy power \cite{Sha}, valid when $p=1$, since it has been constructed to be linear in time in correspondence to the Barenblatt solution of the nonlinear diffusion equation \fer{poro}. Indeed, owing to definition \fer{epo}, it is a simple exercise to show that
 \be\label{lin2}
\mN_p(M_p(t)) = \mN_p(\tilde M_p)\cdot t.
 \ee
Definition \fer{epo} requires $p >(n-2)/n$, in which case $2/n + p-1
>0$.
The range of the parameter $p$  for which we can introduce our
notion of R\'enyi entropy power, coincides with the range for which
there is mass conservation for the solution of \fer{poro}
\cite{BDV}. This range includes the cases in which the Barenblatt
has bounded second moment, since $(n-2)/n < n/(n+2)$.

Definition \fer{epo} allows to extend the concavity property of Costa \cite{Cos}, relative to the linear diffusion equation,  to R\'enyi entropy power of the solution to the nonlinear diffusion equation \fer{poro}, as soon as the second moment of the Barenblatt solution is bounded \cite{ST}. 
The precise result is the following. 
  Let  $p>n/(n+2)$ and let $u(\cdot,t)$ be probability densities in
  $\R^n$ solving the initial value problem for equation \fer{poro} for $t>0$.
  Then the  $p$-th R\'enyi entropy power
  defined by \fer{epo} has the \emph{concavity property} so that
   \be\label{conc-p}
   \frac{d^2}{d
    t^2}\mN_p(v(\cdot,t)) \le 0.
     \ee
Moreover, equality in \fer{conc-p} is achieved if and only if $v(t)$ is the Barenblatt solution \fer{ba-self}. 
Like in Shannon's case, inequality \fer{conc-p} lieds to sharp
isoperimetric inequalities. Under the allowed conditions on the parameter $p$,  the first (nonincreasing) derivative of $p$-th R\'enyi entropy power along the solution to the nonlinear diffusion gives  the scale invariant \emph{isoperimetric inequality} 
  \be\label{b5}
  \mN_p(f) \, \mI_p(f) \ge \mN_p(\tilde
  M_p) \, \mI_p(\tilde M_p) =
    \gamma_{n,p}.
  \ee
Indeed, it is immediate to show that the product in \fer{b5} is invariant under dilation, which allows to reckon explicitly the value of the constants $\gamma_{n,p}$ by using the same argument of Section \ref{sec:heat} \cite{ST}. Clearly, the limit value is obtained in correspondence to the Barenblatt profile \fer{ba}.

Inequality \fer{b5} can be rewritten in a form more suitable to
functio\-nal ana\-ly\-sis. Indeed, when $f(x)$ be a probability density in $\R^n$, and
 $p > n/(n+2)$, the isoperimetric inequality \fer{b5} takes the form
 \be\label{gn}
\int_{\R^n} \frac{|\nabla f^p(x)|^2}{f(x)} \, d x \ge \gamma_{n,p}
\left( \int_{\R^n} f^p(x) \, d x \right)^{\frac{ 2+ 2n
(p-1)}{n(p-1)}}.
 \ee
By setting into \fer{gn}
 \[
 f^{p-1/2} = \frac u{\| u\|_{L^{2q}(\R^n)}}, \quad q = \frac 1{2p-1},
 \]
one obtains that $f$ is a probability density in $\R^n$, and
 \[
\int_{\R^n} \frac{|\nabla f^p(x)|^2}{f(x)} \, d x = 4p^2q^2\, \int_{\R^n} {|\nabla u(x)|^2}\, d x  \left(\int_{\R^n} u^{2q}\, d x\right)^{-1/q}.
 \]
Hence one realizes that inequality \fer{gn}  is equivalent to one of the two following Gagliardo--Nirenberg inequalities \cite{DT2}. 
If $(n-1)/n \le p<1$ then \fer{gn} is equivalent to
 \be\label{GNa}
\| u\|_{L^{2q}(\R^n)} \le K_{GN} \|\nabla u\|_{L^2(\R^n)}^\theta \| u\|_{L^{q+1}(\R^n)}^{1-\theta}, 
 \ee
where
\[
1 < q \le \frac n{n-2}, \quad \theta = \frac nq\frac{q-1}{n +2 -q(n-2)}.
\]
If $p>1$,  then \fer{gn} is equivalent to
 \be\label{GNb}
\| u\|_{L^{q+1}(\R^n)} \le K_{GN} \|\nabla u\|_{L^2(\R^n)}^\theta \| u\|_{L^{2q}(\R^n)}^{1-\theta}, 
 \ee
where
\[
0<q<1, \quad  \theta = \frac nq\frac{q-1}{n +2 -q(n-2)}.
\]
In particular,
\emph{Sobolev} inequality \cite{Aub,Tal} is obtained when $\theta =1$, $p = 1/n$ and $2q = n/(n-2)$.  Hence, Sobolev inequality with the
sharp constant is a consequence of the concavity of  R\'enyi entropy
power of parameter $ p = (n-1)/n$, when $n >2$, and represents the threshold case for the validity of the method. These
Gagliardo-Nirenberg type inequalities  with sharp constants, have been first obtained by Del Pino and Dolbeault \cite{DD}, and later on by
Cordero-Erausquin, Nazaret, and Villani, \cite{CNV} with different
methods,  without noticing their connection with classical concepts of information theory. 
However,  in the one-dimensional case, these connections have been outlined in \cite{LYZ}, where a more general class of generalized Fisher-type measures of information have been considered, together with the corresponding isoperimetric inequalities. 

\section{Conclusions}

In this note we presented  the principal connections between information-theoretic inequalities for entropies and interpolation inequalities of Gagliardo--Nirenberg type. It is interesting to remark that the proof of most information-type inequalities takes advantage of the time-evolution of entropies along the solutions to diffusion equations, which are known to converge towards self-similar ones. The main example is furnished by the result of concavity presented in Section \ref{sec2}, where the classical sharp Sobolev inequality is derived as a particular case of the concavity property of R\'enyi entropy power of order $1-1/n$.  Like in the linear case, the knowledge of isoperimetric inequalities for entropies is a powerful tool that, among other applications, allows to obtain sharp convergence results on the large-time behavior of the solution to nonlinear diffusion equations \cite{CTo}.  

\vskip 0,5cm
\section*{Acknowledgement} This work has been written within the
activities of GNFM group  of INdAM (National Institute of
High Mathematics), and partially supported by  MIUR project ``Optimal mass
transportation, geometrical and functional inequalities with applications''.
The research was partially supported by the Italian Ministry of Education, University and Research (MIUR): 
Dipartimenti di Eccellenza Program (2018--2022) - Dept. of Mathematics ``F. Casorati'', University of Pavia. 

I wish to thank the editors  of this volume for inviting me to contribute to the memory of the notable figure of Emilio Gagliardo, who was an esteemed colleague in my Department for more than twenty years, starting in November 1975. 


\begin{thebibliography}{}
\bibitem{Ag1}
M. Agueh, N. Ghoussoub, X. Kang, Geometric inequalities via a general comparison principle for interacting gases. \emph{Geom. Funct. Anal.}, \textbf{14}
(2004) 215--244

\bibitem{Agu}
M. Agueh, Gagliardo--Nirenberg inequalities involving the gradient $L^2$-norm. \emph{C. R. Acad. Sci. Paris, Ser. I}, \textbf{346}  (2008) 757--762 

\bibitem{Aub}
T. Aubin, Probl\'emes isop\'erim\'etriques et espaces de Sobolev. \emph{J. Differential Geometry}, \textbf{11} (1976) 

\bibitem{Ba1}
G.I. Barenblatt,  On some unsteady motions of a liquid and gas in a
porous medium. \emph{Akad. Nauk SSSR. Prikl. Mat. Meh.}, \textbf{16}
(1952) 67--78.


\bibitem{Bar} G.I. Barenblatt, \emph{Scaling, self-similarity and intermediate asymptotics}.
Cambridge Univ. Press, Cambridge 1996.

\bibitem{B2}
J. Bennett, {Heat-flow monotonicity related to some inequalities in
euclidean analysis}. \emph{Contemporary Mathematics}, \textbf{505}  (2010),
85--96.

\bibitem{BB}
J. Bennett and N. Bez, {Closure properties of solutions to heat
inequalities}. \emph{J. Geom. Anal.}, \textbf{19}  (2009),  584--600.

\bibitem{BBC}
J. Bennett, N. Bez and A. Carbery, {Heat-flow monotonicity
related to the Hausdorff--Young inequality}, \emph{Bull. London Math. Soc.}
\textbf{41} (6) (2009), 971--979.

\bibitem{B1}
J. Bennett, A. Carbery, M. Christ and T. Tao, {The
Brascamp--Lieb inequalities: finiteness, structure and extremals}.
\emph{Geometric And Functional Analysis}, \textbf{17}  (2007), 1343--1415.


 \bibitem{Bla}
N.M. Blachman, The convolution inequality for entropy powers.
\emph{IEEE Trans. Inform. Theory}, \textbf{11} (1965) 267--271.

\bibitem{BDV}
A. Blanchet, M. Bonforte, J. Dolbeault, G. Grillo, and J. L V\'azquez,
Asymptotics of the fast diffusion equation via entropy estimates.
\emph{Arch. Ration. Mech. Anal.},  \textbf{191} (2009) 347--385.

\bibitem{BDGV}
M.~{B}onforte, J.~{D}olbeault, G.~{G}rillo, and J.-L. {V}{\'a}zquez,
  {S}harp rates of decay of solutions to the nonlinear fast diffusion equation
  via functional inequalities. \emph{Proc. Natl. Acad. Sci. USA}, \textbf{107} (2010)
  16459--16464.
  
\bibitem{CF}
E.A. Carlen, and A. Figalli,
Stability for a GNS inequality and the Log-HLS inequality, with application to the critical mass Keller--Segel equation. \emph{Duke Math. J.},
\textbf{162} (3) (2013) 579--625.
  

\bibitem{CaTo} J. A. Carrillo, and G. Toscani, \newblock Asymptotic $L^1$-decay of
solutions of the porous medium equation to self-similarity.  {\em Indiana Univ. Math. J.},  {\bf 49} (2000)
113--141.  

\bibitem{CTo}
J.A. Carrillo, and G. Toscani, R\'enyi entropy and improved rate of convergence to self-similarity for nonlinear diffusion equations. \emph{Nonlinearity}, \textbf{27} (2014)  3159--3177.

\bibitem{CNV}
D. Cordero-Erausquin, B. Nazaret, and C. Villani,  A mass--transportation
approach to sharp Sobolev and Gagliardo--Nirenberg inequalities.
\emph{Advances in Mathematics}, \textbf{182}  (2004) 307--332.

\bibitem{Cos}
 M. Costa, A new entropy power inequality.   \emph{IEEE Trans. Inform. Theory},
\textbf{31} (1985) 751--760.

\bibitem{CHV}
J. Costa, A. Hero, and C. Vignat, On solutions to multivariate
maximum alpha-entropy problems. \emph{Lecture Notes in Computer
Science}, \textbf{2683}, no. EMMCVPR 2003, Lisbon, 7-9 July 2003, (2003)
211--228.

 \bibitem{CT}
T.M. Cover,  and J.A. Thomas, \emph{Elements of information theory}.
 Wiley \& Sons, New York,
1992.



\bibitem{DD}
M. Del Pino, and J. Dolbeault, Best constants for Gagliardo--Nirenberg
inequalities and application to nonlinear diffusions. \emph{J. Math. Pures
Appl.},  \textbf{81} (2002) 847--875.

\bibitem{Dem}
A. Dembo,  A simple proof of the concavity of the entropy power with
respect to the variance of additive normal noise. \emph{IEEE Trans.
Inform. Theory}, \textbf{35} (1989)  887--888.


 \bibitem{DCT}
  A. Dembo, T. Cover, and J. Thomas, Information theoretic
inequalities. \emph{IEEE Trans. Inform. Theory},  \textbf{37} (1991)
1501--1518.

\bibitem{DT1}
 J. Dolbeault, and G. Toscani, Best matching Barenblatt profiles are delayed. \emph{J. Phys. A: Math. Theor.}, \textbf{48}  (2015) 065206.

\bibitem{DT3}
J. Dolbeault, and G. Toscani, Stability results for logarithmic Sobolev and Gagliardo--Nirenberg inequalities. \emph{Int. Math. Res. Notices}, rnv 131 (2015). 

\bibitem{DT2}
 J. Dolbeault, and G. Toscani, Nonlinear diffusions: extremal properties of Barenblatt profiles, best matching and delays.  \emph{Nonlinear Anal. Series A}, \textbf{138} (2016) 31--43.
 
\bibitem{Eva}
L.C. Evans, \emph{Partial differential equations}. Graduate Studies in Mathematics 19. AMS, Providense, RI 2010.

\bibitem{Gagl} 
E. Gagliardo, Propriet$\grave{\rm a}$ di alcune classi di funzioni in pi$\grave{\rm u}$ variabili. \emph{Ricerche Mat.}, \textbf{7} (1958) 102--137.

\bibitem{GSV}
D. Guo, S. Shamai, and S. Verd\'u, Mutual information and minimum
mean-square error in Gaussian channels. \emph{IEEE Trans. Inf.
Theory}, \textbf{51}, (4) (2005) 1261--1282.

\bibitem{GSV2}
D. Guo, S. Shamai, and S. Verd\'u, A simple proof of the
entropy-power inequality . \emph{IEEE Trans. Inf. Theory},
\textbf{52}, (5) (2006) 2165--2166.

\bibitem{McK}
\textsc{H.P. McKean Jr.},  {Speed of approach to equilibrium
for Kac's caricature of a Maxwellian gas}. \emph{Arch. Rat. Mech. Anal.},
\textbf{21}  (1966) 343--367.

\bibitem{Lev}
 H.A. Levine, An estimate for the best constant in a Sobolev inequality involving three integral norms.  \emph{Ann. Mat. Pura Appl.},  \textbf{124} (1980) 181--197.

\bibitem{Lin}
Yu.V. Linnik, {An information-theoretic proof of the central limit
theorem with the Lindeberg condition}. \emph{Theory Probab. Appl.},
\textbf{4}  (1959)  288-299.

\bibitem{LT}
P.L.Lions, and G.Toscani, A sthrenghtened central limit theorem for smooth densities. \emph{J. Funct. Anal.}, \textbf{128} (1995) 148--167.

 \bibitem{LYZ}
E. Lutwak, D. Yang, and G.  Zhang,
 Cram\'er-Rao and moment-entropy inequalities for R\'enyi entropy and generalized Fisher
 information. \emph{IEEE Trans. Inform. Theory}, \textbf{51} (2005)
473--478.

\bibitem{LYZ2}
E. Lutwak, D. Yang, and G. Zhang, Moment-entropy inequalities for a
random vector. \emph{IEEE Trans. Inform. Theory}, \textbf{53} (2007)
1603--1607.

\bibitem{Nash}
J.F. Nash, Continuity of solutions of parabolic and elliptic equations. \emph{Amer. J.
Math.},  \textbf{80} (1958) 931--954.

\bibitem{Nire}
L. Nirenberg, On elliptic partial differential equations. \emph{Ann. Scuola Norm. Pisa}, \textbf{13} (1959) 116--162.



\bibitem{Ren}
A. R\'enyi,  On measures of entropy and information. \emph{Proc. Fourth Berkeley Symp.
Math. Statist. Prob. 1}, 547--561. University of California Press, Berkeley  1961.

\bibitem{Rio}
O. Rioul, Information theoretic proofs of entropy power inequalities. \emph{IEEE
Trans. Inform. Theory}, \textbf{57} (1) (2011)  33--55.

\bibitem{ST}  G. Savar\'e, and G. Toscani,  The concavity of R\'enyi entropy power. \emph{IEEE Trans. Inform. Theory}, \textbf{60} (2014) 2687--2693. 


\bibitem{Sha}
C.E. Shannon, {A mathematical theory of communication}. \emph{Bell
Syst. Tech. J.}, \textbf{27} Jul.(1948), 379--423; Oct. (1948)
623--656.

\bibitem{Sta}
 A.J. Stam, {Some inequalities satisfied by the quantities of
information of Fisher and Shannon}. \emph{Inf. Contr.}, \textbf{2}
(1959) 101--112.

\bibitem{Tal}
G. Talenti, Best constant in Sobolev inequality. \emph{Ann. Mat. Pura
Appl.}, \textbf{110} (1976) 353--372.

\bibitem{Tos2}
G. Toscani,  {An information-theoretic proof of Nash's inequality}.
\emph{Rend. Lincei Mat. Appl.},  \textbf{24} (2013) 83--93.

\bibitem{Tos3}
 G. Toscani, {Lyapunov
functionals for the heat equation and sharp inequalities}. \emph{Atti
Acc. Peloritana Pericolanti, Classe Sc. Fis. Mat. e Nat.},
\textbf{91} (2013)  1--10.

\bibitem{Tos4}
 G. Toscani, Heat equation and convolution inequalities.  \emph{Milan J. Math.}, \textbf{82} (2014) 183--212. 

\bibitem{Tos5}
 G. Toscani, R\'enyi entropies and nonlinear diffusion equations. \emph{Acta. Appl. Math.}, \textbf{132} (2014) 595--604.

\bibitem{Vaz}
J.L. Vazquez, \emph{The porous medium equation: Mathematical
theory}. Oxford University Press, Oxford 2007.

 \bibitem{Vil}
C. Villani, A short proof of the concavity of entropy power. \emph{IEEE Trans. Inform. Theory}, 
\textbf{46} (2000) 1695--1696.

\bibitem{ZK} Ya.B. Zeldovich, and A.S. Kompaneetz,  Towards a theory of
heat conduction with thermal conductivity depending on the
temperature. Collection of Papers Dedicated to 70th Birthday of
Academician A. F. Ioffe, \emph{Izd. Akad. Nauk SSSR}, Moscow, (1950)
61--71.

\end{thebibliography}


\end{document}